\documentclass[a4paper,12pt]{article}

\usepackage{amsmath,amssymb,latexsym,graphics,epsfig}
\usepackage{hyperref}
\usepackage{color}
\usepackage{amsthm}
\usepackage{graphicx,url}
\usepackage{amsopn}
\usepackage{amssymb}
\usepackage[english]{babel}
\usepackage{soul}
\usepackage{tablefootnote}
\usepackage[normalem]{ulem}

\definecolor{Dominika}{RGB}{255,0,0}

\newtheorem{Theorem}{Theorem}[section]
\newtheorem{Lemma}[Theorem]{Lemma}

\newtheorem{Proposition}[Theorem]{Proposition}

\newtheorem{Example}[Theorem]{Example}

\DeclareMathOperator{\dist}{dist}

\DeclareMathOperator{\ecc}{ecc}

\def\Z{\mathbb Z}

\def\j{\mbox{\boldmath $j$}}

\def\s{\mbox{\boldmath $s$}}

\def\vec0{\mbox{\boldmath $0$}}

\def\A{\mbox{\boldmath $A$}}
\def\B{\mbox{\boldmath $B$}}

\def\I{\mbox{\boldmath $I$}}

\def\Z{\mathbb Z}

\begin{document}

 \title{Some inner metric parameters of a digraph:
 Iterated line digraphs and integer sequences
	}
	\author{N. H. Bong$^a$, C. Dalf\'o$^b$,  M. A. Fiol$^c$, and D. Z\'avack\'a$^{d}$\\
		\\
  {\small $^a$Dept. of Mathematics, Delaware University, Newark, Delaware, USA} \\ 	
		{\small {\tt{nhbong@udel.edu}}}\\
{\small $^b$Dept. de Matem\`atica, Universitat de Lleida, Igualada (Barcelona), Catalonia}\\
        {\small {\tt cristina.dalfo@udl.cat}}\\
{\small $^{c}$Dept. de Matem\`atiques, Universitat Polit\`ecnica de Catalunya, Barcelona, Catalonia} \\
{\small Barcelona Graduate School of Mathematics} \\
{\small  Institut de Matem\`atiques de la UPC-BarcelonaTech (IMTech)}\\
		{\small {\tt miguel.angel.fiol@upc.edu} }\\
{\small $^d$Dept. of Applied Informatics, Comenius University, 
  Bratislava, Slovakia} \\
		{\small {\tt  dominika.mihalova@fmph.uniba.sk }}
	}

\date{}

\maketitle

\begin{abstract}
In this paper, we first give a new result characterizing the strongly connected digraphs with a diameter equal to that of their line
digraphs. 
 Then, we introduce the concepts of the inner diameter and inner radius of a digraph and study
their behaviors in its iterated line digraphs.
Furthermore, we provide a method to characterize sequences of integers (corresponding to the inner diameter or the number of vertices of a digraph and its iterated line digraphs) that satisfy some conditions. Among other examples, we apply the method to the cyclic Kautz digraphs,
square-free digraphs, and the subdigraphs of De Bruijn digraphs. Finally, we present some tables with new sequences that do not belong to The On-Line Encyclopedia of Integer Sequences.
\end{abstract}




\section{Preliminaries}

Let us first introduce some basic notation and results.
A digraph $G=(V,E)$ consists of a (finite) set
$V=V(G)$ of vertices and a multiset $E=E(G)$ of arcs (directed edges) between vertices
of $G$. Here, we also assume that $V$ and $E$ are finite set and multiset, respectively. As the initial and final vertices of an arc are not necessarily different, the
digraphs may have \emph{loops} (arcs from a vertex to itself) 
and \emph{multiple arcs}, that is, there can be more than one arc from one vertex
to another. If $a=(u,v)$ is an arc from $u$ to  $v$, then vertex $u$ (and arc $a$)
is {\em adjacent to} vertex $v$, and vertex $v$ (and arc $a$) is {\em adjacent from} $u$. 
The {\em converse digraph} $\overline{G}$
is obtained from $G$ by reversing the orientation of all its arcs.
Let $G^+(v)$ and $G^-(v)$ denote the sets of arcs adjacent from and to vertex
$v$, respectively. Their cardinalities are the {\em out-degree} $\delta^+(v)=|G^+(v)|$ and {\em in-degree} $\delta^-(v)=|G^-(v)|$ of $v$. A digraph $G$ is $\delta$\emph{-regular} if $\delta^+(v)=\delta^-(v)=\delta$ for all $v\in V$.

In the line digraph $LG$ of a digraph $G$, each vertex of $LG$ represents an
arc of $G$, that is, $V(LG)=\{uv : (u,v)\in E(G)\}$; and vertices $uv$ and $wz$ of $L(G)$ are adjacent if and only if $v=w$, namely, when arc $(u,v)$ is adjacent to arc $(w,z)$ in $G$.
The $k$-iterated line digraph $L^kG$ is recursively defined as $L^0G=G$ and $L^kG=LL^{k-1}G$ for $k\geq 0$.
A digraph $G$ is \emph{periodic} if $L^mG=L^{m+k}G$ for some integers $m$ and $k$, with $k>0$, and the smallest value of $k$ is called the \emph{period} of $G$. Moreover, it was determined when, for two digraphs $G_1$ and $G_2$, there exist integers $m$ and $n$ such that $L^mG_1=L^nG_2$. For more information on periodic digraphs, see Hemminger \cite{He74}.

It can easily be seen that every vertex of $L^kG$ corresponds to a (directed) walk
$v_0,v_1,\ldots,$ $v_k$ of length $k$ in $G$, where $(v_{i-1},v_{i})\in E$ for $i=1,\ldots,k$.
If there is at most one arc between pairs of vertices and $\A$ is the adjacency matrix of $G$, then the $uv$-entry of the power $\A^k$, denoted by $a_{uv}^{(k)}$, corresponds to the number of $k$-walks from vertex $u$ to vertex $v$ in $G$. Furthermore, the order $n_k$ of $L^kG$ turns out to be
\begin{equation}\label{orderL^kG}
n_k=\j^{\top}\A^k\j,
\end{equation}
where $\j$ stands for the all-$1$ \textcolor{red}{(column)} vector. If multiple arcs exist between pairs of vertices, then the corresponding entry in the matrix is not 1, but the number of these arcs.
If $G$ is a $d$-regular digraph with $n$ vertices, then its line
digraph $L^kG$ is
$d$-regular with $n_k=d^kn$ vertices.

The distance from vertex $u$ to vertex $v$ in $G$, denoted by $\dist_G(u,v)$, is the length of a shortest path from $u$ to $v$. Notice that, in general, this does not define a metric since it is possible that  $\dist_G(u,v)\neq \dist_G(v,u)$.
Recall also that a digraph $G$ is \emph{strongly connected} if there is a
(directed) walk between every
pair of its vertices $u$ and $v$, that is, $\dist_G(u,v)<\infty$. Note that a digraph is strongly connected if and only if its line digraph is also connected.
Generalizing this concept, we say that $G$ is {\em unilaterally connected} if, for any pair of vertices $u$ and $v$, either $\dist_G(u,v)<\infty$ or $\dist_G(v,u)<\infty$, see Dalf\'o and Fiol \cite{df14}.
If the digraph $G$ is not strongly connected but its underlying graph $UG$ is connected, then $G$
is called \emph{weakly
connected}. 

For the concepts and results on digraphs not presented here, 
see, for instance, Bang-Jensen and Gutin \cite{BJGu09}, 
Chartrand and Lesniak~\cite{cl96}, or Diestel~\cite{d10}. 


This paper is organized as follows. In Section~\ref{sec:stand-diam}, we recall some known facts on standard metric parameters, and we give a new result characterizing the strongly connected digraphs with a diameter equal to that of their line digraphs.
In Section~\ref{sec:inner-diam}, we prove one of our main results and mention some possible examples. 
Furthermore, 
in Section 
\ref{sec:sequences}, we characterize sequences of integers that satisfy some conditions. These correspond to the numbers of vertices of iterated line digraph, such as the cyclic Kautz digraphs, the unicyclic digraphs, and the acyclic digraphs, considered in Subsections \ref{sec:ex-CK}, \ref{sec:ex-unicyclic}, and \ref{sec:ex-acyclic}, respectively.
Finally, we present some tables with new sequences that do not belong to The On-Line Encyclopedia of Integer Sequences.

\section{Metric parameters}
\label{sec:metric-param}

\subsection{Standard metric parameters}
\label{sec:stand-diam}

If $G$ is a strongly connected digraph, different from a directed cycle, with
diameter $D$, then its line digraph $L^kG$ has diameter $D+k$. See Fiol, Yebra, and Alegre~\cite{fya84} for more details.
The interest of the line digraph technique is that it allows us to obtain digraphs with small diameters and large connectivity.
To compare the line digraph technique and other techniques to get digraphs with minimum diameter, see Miller, Slamin, Ryan, and Baskoro~\cite{MiSlRyBa13}. Since these techniques are related to the degree/diameter problem, we also refer to the comprehensive survey on this problem by Miller and \v{S}ir\'{a}\v{n}~\cite{ms13}.

Note that the directed cycle (also called {dicycle}) and its line digraph are isomorphic, and so the diameters of both digraphs coincide. In fact, we show that the converse is also true, which, as far as we know, is a new result.

\begin{Proposition}
	\label{lemma:digrafs}
	Let $G$ be a strongly connected digraph. Then, $D(LG)=D(G)$ if and only if $G$ is a directed cycle.
\end{Proposition}

\begin{proof}
If $G$ is a directed cycle, then $LG\cong G$ and their diameters coincide. Conversely, let us prove first that if the diameters coincide, $G$ must contain a directed cycle of length $D+1$.
Consider two vertices $u$ and $v$ at maximum distance $D$ with shortest path $u_0(=u),u_1,\ldots,u_D(=v)$. Since $G$ is strongly connected, there must be an arc going to $u$, called $(u',u)$, and an arc going from $v$ called $(v,v')$. Notice that neither $(u',u)$ nor $(v,v')$ can belong to the shortest path from $u$ to $v$; otherwise, the distance between $u$ and $v$ in $G$ would be smaller than $D$. Moreover, if we had that $(u',u)\neq (v,v')$, then the diameter of the line digraph $LG$ would be $D+1$ because of the shortest path  $u'u_0,u_0u_1,\ldots,u_Dv'$, which is a contradiction with the hypothesis. Consequently, $(u',u)=(v,v')$, forming the claimed directed cycle.
	
Furthermore, we show that $G$ can only be a directed cycle. Suppose a directed cycle is a (proper) subdigraph of $G$. Then,  there must be either an arc going to or (an arc going) from a vertex $u_i$ in the directed cycle. Without loss of generality, we can assume the first case, with the incoming arc $(w_i, u_i)$  not belonging to the cycle. (Otherwise, we simply reason with the converse digraph $\overline{G}$, which would be a cycle if and only if  $G$ is.) Now, since $G$ is strongly connected, so is the line digraph $LG$. Hence, there should be a directed path from vertex $u_i$ to vertex $u_{i-1}$ in order to have $\dist_{LG}(w_iu_i,u_{i-1}u_{i})=\dist_G(u_i,u_{i-1})+1\le D$. Thus, $\dist_G(u_i,u_{i-1})<D$, and so there is a shortest path from $u_i$ to $u_{i-1}$ of length smaller than $D$. Let $u_j,u_{j+1},\ldots,u_{i-1}$ be the last vertices of that shortest path belonging to the directed cycle, and let $(w_j, u_j)$ be the arc going to $u_j$. Note that $u_j\neq u_{i+1}$. Now, we consider the following cases (see Figure \ref{fig:demo-digrafs}):

\begin{figure}[!ht]
	\begin{center}
		\includegraphics[width=12cm]{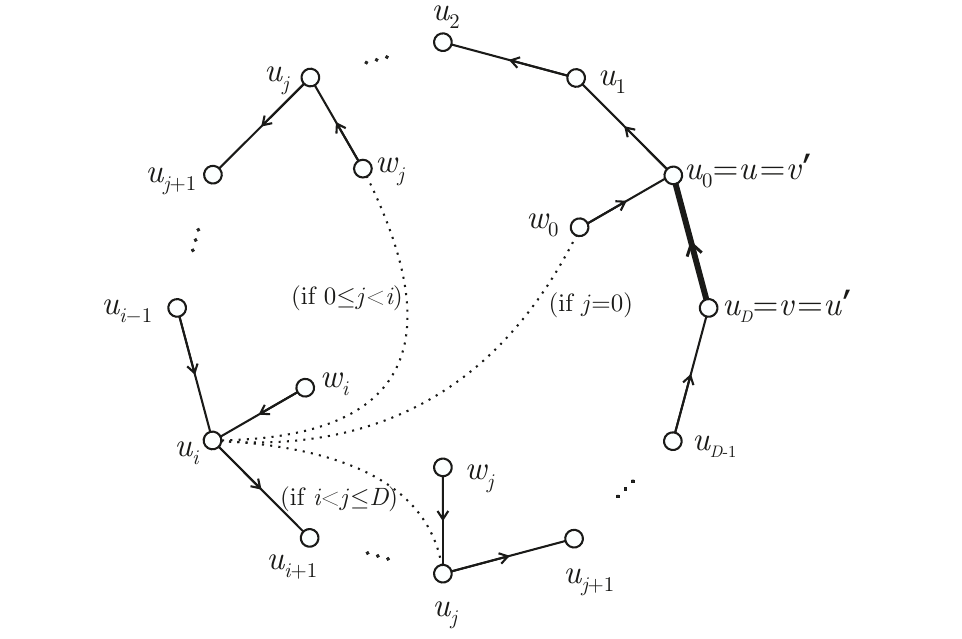}
	\end{center}
	\vskip-.5cm
	\caption{Scheme of the proof of Proposition \ref{lemma:digrafs}.}
	\label{fig:demo-digrafs}
\end{figure}

\begin{itemize}
\item[$(i)$]
If $i< j\leq D$, the shortest path $u_i\rightarrow u_{i-1}$ would contain the arc $(v,u)$ and, hence, there would be a directed path from $u$ to $v$ of length smaller than $D$, a contradiction.
\item[$(ii)$] If $j=0$, then we would have the arc $(w_0,u_0)=(w_0,u)\neq (v,u)$ contradicting that there is only one arc going to $u$. Recall that if we had that $(w_0, u)\neq (v,u)$, then the diameter of the line digraph $LG$ would be $D+1$.
\item[$(iii)$]
Thus, it should be $0< j <i$, but  then we can use the same reasoning $(i)$ and $(ii)$ with $u_j$ instead of $u_i$ to show that, in the  shortest directed path from $u_j$ to $u_{j-1}$, there exists a vertex $u_k$, with     $0<k<j$.
\item[$(iv)$] Repeating these steps, we eventually would reach a vertex $u_l$ such that the path from $u_l$ to $u_{l-1}$ would go through a vertex $u_h$, with only the forbidden cases $(i)$ $l<h\leq D$, or $(ii)$ $l=0$, a contradiction.
\end{itemize}
This completes the proof.
\end{proof}

\subsection{Inner metric parameters}
\label{sec:inner-diam}

Let $G=(V,E)$ be a (not necessarily strongly connected) digraph. 
We define the \emph{inner distance} between vertices $u$ and $v$ as the (standard) distance whenever it is not infinite. So, the inner distance is not defined for pairs of vertices with no path between them. From now on, we are going to refer to the inner distance (or distance) when it is not infinite. 

\begin{figure}[!ht]
	\begin{center}
		\includegraphics[width=4cm]{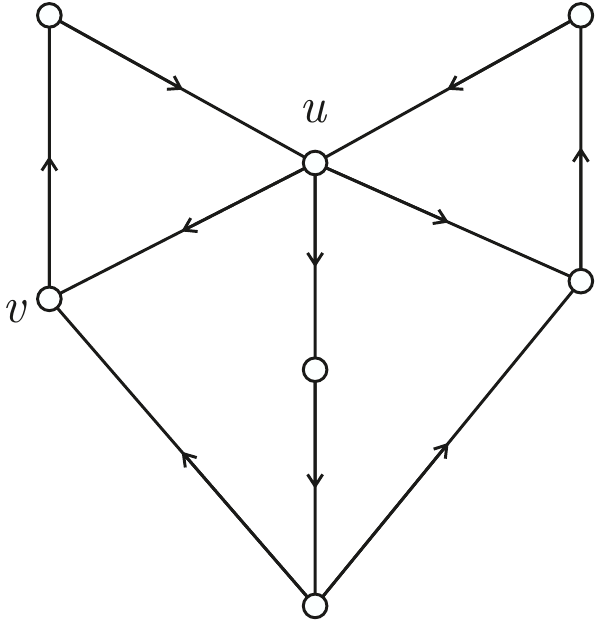}
	\end{center}
 \vskip-.5cm
	\caption{A digraph $G$ with different in- and out-radius.}
	\label{fig:radii}
\end{figure}

In consequence, we now
define the {\em inner out-eccentricity} of a vertex $u\in V$ as
$\ecc^+(u)=\max_{v\in V}\{\dist_G(u,v): \dist_G(u,v)<\infty \}$. Similarly, the {\em inner in-eccentricity} of  $u\in V$ is
$\ecc^-(u)=\max_{v\in V}\{\dist_G(v,u): \dist_G(v,u)<\infty \}$. From this, we can define the following parameters:
\begin{itemize}
\item
{\em Inner in-radius}: $r^-(G)=\displaystyle\min_{u\in V}\ecc^-(u)$.
\item
{\em Inner (out-)radius}: $r(G)=r^+(G)=\displaystyle\min_{u\in V}\ecc^+(u)$.
\item
{\em Inner in-diameter}: $d^-(G)=\displaystyle\max_{u\in V}\ecc^-(u)$.
\item
{\em Inner (out-)diameter}: $d(G)=d^+(G)=\displaystyle\max_{u\in V}\ecc^+(u)$\\ 
\mbox{\hskip 5cm} $=\max_{u,v\in V}\{\dist_G(u,v): \dist_G(u,v)<\infty \}$.
\item
{\em Inner mean distance}:
$\displaystyle\overline{d}=\frac{1}{|U|}\sum_{(u,v)\in U}\dist(u,v),$\\
where $U=\{(u,v):\dist(u,v)<\infty\}$.
\end{itemize}
From these definitions, it is clear that $r^-(G)=r^+(\overline{G})$ and 
$d^-(G)=d^+(\overline{G})$. 
As we show in the next lemma, the inner in- and out-diameter coincide. This is not the case for the inner in- and out-radius. For instance, if $G$ has some vertex $u$ with in-degree 0 but no vertex with out-degree 0, then $r^-(G)=\ecc^-(u)=0$, but $r^+(G)>0$. 
This can occur even if $G$ is strongly connected. For instance, consider the digraph with vertices $u,v,w,x$ and arcs $uv$, $vw$, $wx$, $vu$, $wu$, $xu$. In this case, $r^+(G)=\ecc^+(v)=2$ but $r^-(G)=\ecc^-(u)=1$. Another example is the digraph $G$ of Figure \ref{fig:radii}, where $r^+(G)=\ecc^+(u)=2$ but $r^-(G)=\ecc^-(v)=3$.
Further, in this section, we prove a stronger result showing the existence of strongly connected digraphs for any values of inner in- and out-radius.

\begin{Lemma}
\label{lem:in-out-d}
For any digraph $G$, the inner in- and out-diameters coincide: $d^-(G)=d^+(G)=d(G)$.    
\end{Lemma}
\begin{proof}
Let $u,v$ be two vertices at maximum distance $\dist(u,v)=d^+(G)=\ell$, so that we have the longest path $u_0(=u),u_1,\ldots,u_{\ell}(=v)$. This implies that $\ecc^-(v)\ge \ell$ and, hence,
$d^-(G)\ge \ecc^-(v)\ge d^+(G)$. Analogously, if we consider two vertices $x,y$ such that $\dist(y,x)=d^-(G)$, we obtain that $d^+(G)\ge d^-(G)$, and the result follows.
\end{proof}

\begin{figure}[!ht]
	\begin{center}
		\includegraphics[width=16cm]{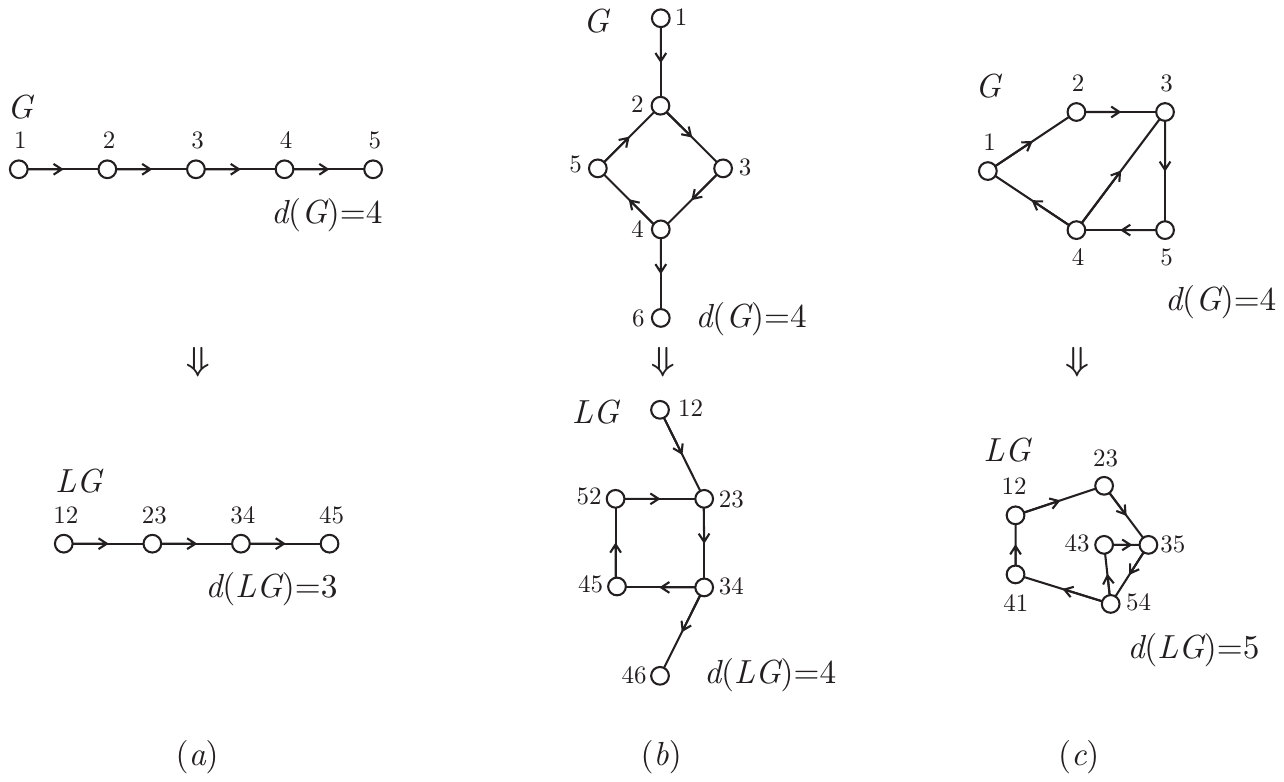}
	\end{center}
	\vskip-4.25cm
	\caption{Examples of cases $(a)$, $(b)$, and $(c)$ in the proof of Proposition \ref{prop:3diam}.}
	\label{fig:ex-inner-diam}
\end{figure}

As mentioned before, for any strongly connected digraph $G$ different from a directed cycle, $D(LG)=D(G)+1$. However, this is not necessarily true for the inner diameter when $G$ is not strongly connected. See, for example, the  digraph $G$ represented in Figure \ref{fig:ex-inner-diam}$(a)$ has inner diameter $d(G)=4$, and its line digraph $LG$ has inner diameter $d(LG)=3$. In general, we have the following result.




\begin{Proposition}
\label{prop:3diam}
Let $G$ be a digraph with inner diameter $d(G)$. Then, the inner diameter $d(LG)$ of its line digraph $LG$ satisfies
\begin{equation}
d(G)-1\leq d(LG)\leq d(G)+1.
\label{d(LG)}
\end{equation}
\end{Proposition}

\begin{proof}
	Let us consider a shortest path in $G$ between two vertices $u$ and $v$ at distance $\ell$:
	$$
	u(=u_0), u_1,\ldots,v(=u_{\ell}).
	$$
	Then, there are three different cases: An example of these three cases is illustrated in Figure \ref{fig:ex-inner-diam}.)
	\begin{enumerate}
		\item[$(i)$] If there is neither a vertex adjacent to $u$, nor a vertex adjacent  from $v$ (that is, $G^-(u)=G^+(v)=\emptyset$), the induced shortest path in $LG$
		$$
		u_0u_1,u_1u_2,\ldots,u_{{\ell}-1}u_{\ell}
		$$
		implies that $\dist_{LG}(u_0u_1,u_{{\ell}-1}u_{\ell})={\ell}-1$.
		\item[$(ii)$] If there exists a vertex adjacent to $u$ (that is, $u_{-1}\in G^-(u)$), but there is no vertex adjacent from $v$ (that is, $G^+(v)=\emptyset$), or vice versa, if there is no vertex adjacent to $u$ (that is, $G^-(u)=\emptyset$), but there exists a vertex adjacent from $v$  (that is, $u_{{\ell}+1}\in G^+(v)$), then the corresponding induced shortest paths in $LG$ are
		$$
		u_{-1}u_0,u_0u_1,u_1u_2\ldots,u_{{\ell}-1}u_{\ell} \quad \mbox{and} \quad u_0u_1,u_1u_2\ldots,u_{{\ell}-1}u_{\ell},u_{{\ell}}u_{{\ell}+1},
		$$
		and indicate that $\dist_{LG}(u_{-1}u_0,u_{{\ell}-1}u_{\ell})=\dist_{LG}(u_0u_1,u_{\ell}u_{{\ell}+1})={\ell}$.
		\item[$(iii)$] If there exist both a vertex $u_{-1}$ adjacent to $u$ and a vertex $u_{{\ell}+1}$ adjacent from $v$, then the induced shortest path in $LG$
			$$
		u_{-1}u_0,u_0u_1,u_1u_2\ldots,u_{{\ell}-1}u_{\ell},u_{{\ell}}u_{{\ell}+1}
		$$
		implies that $\dist_{L(G)}(u_{-1}u_0,u_{{\ell}}u_{{\ell}+1})={\ell}+1$, unless $u_{-1}u_0=u_{{\ell}}u_{{\ell}+1}$ (that is, $u_{-1}=u_{\ell}$ and $u_{0}=u_{{\ell}+1}$), so that $u_0, u_1,\ldots,u_{\ell},u_{{\ell}+1}$ forms a directed cycle in $G$.
		\end{enumerate}
	Consequently,
\begin{enumerate}
\item[$(a)$]
If all the shortest paths of length ${\ell}$ in $G$ are of type $(i)$, then
$$
d(LG)=d(G)-1.
$$
\item[$(b)$]
If there is a shortest path of length ${\ell}$ in $G$ of type $(ii)$, but none of type $(iii)$ (with exception of a directed cycle), 
then
$$
d(LG)=d(G).
$$
\item[$(c)$]
If there is a shortest path of length ${\ell}$ in $G$ of type $(iii)$, different from a directed cycle, 
then
$$
d(LG)=d(G)+1.
$$
\end{enumerate}
\end{proof}

Notice that if $G$ is strongly connected and different from a cycle, then $d(LG)=d(G)+1$ and, in general, $d(L^kG)=d(G)+k$ for $k\ge 0$.

The same reasoning proves that the inner (out-)radius $r^+(LG)$ and inner (in-)radius $r^-(LG)$ satisfy an inequality similar to \eqref{d(LG)}.
\begin{Proposition}
\label{prop:3diam-bis}
Let $G$ be a digraph with inner (out-)radius $r(G)=r^+(G)$ and inner in-radius $r^-(G)$. Then, the corresponding parameters of its line digraph $LG$ satisfies
\begin{align}
r^+(G)-1 & \leq r^+(LG)\leq r^+(G)+1,\\
r^-(G)-1 & \leq r^-(LG)\leq r^-(G)+1.
\label{r(LG)}
\end{align}
In particular, if $G$ is strongly connected and different from a directed cycle, then $r^+(LG)=r^+(G)+1$ and $r^-(LG)=r^-(G)+1$.
\end{Proposition}
 From this result, and in contrast to Lemma \ref{lem:in-out-d}, we now show that  $r^+(G)$ and $r^-(G)$ can take any values.

\begin{Proposition}
For any pair of  positive integers $r_1,r_2$, there exists a  strongly connected digraph $G$ such that $r^+(G)=r_1$  and $r^-(G)=r_2$.
\end{Proposition}
\begin{proof}
Without loss of generality, assume that $r_1\le r_2$ (otherwise, consider the converse $\overline{G}$). Let us first consider the digraph $C_n^*$ with set of vertices labeled in $\Z_n$ and arcs $i\rightarrow i+1 \mod n$
(forming a directed cycle), and $0\rightarrow j$ with $j=2,\ldots,n-1$ (see Figure \ref{fig:Cn*} for the case $n=9$). Then, since $\ecc^+(0)=1$ and $\ecc^-(i)\ge \lceil (n-1)/2\rceil$, we have that $r^+(G)=1$ and $r^-(G)=(n-1)/2$ if $n$ is odd, and $r^-(G)=n/2$ if $n$ is even (in Figure \ref{fig:Cn*}, the label `$i\!\!:\!\!j$'
indicates that $\ecc^-(i)=j$). Then, since $C_n^*$ is strongly connected, its $k$-iterated line digraph $L^kC_n^*$ has inner radii $r^+(G)=1+k$ and 
 $r^-(G)=\lceil (n-1)/2\rceil+k$. Consequently, by taking the values $k=r_1-1$ and $n=2(r_2-r_1)+3$ ($n$ odd) or $n=2(r_2-r_1)+2$ ($n$ even),
 it turns out that $r^+(G)=r_1$ and $r^{\textcolor{red}{-}}(G)=r_2$, as required.
\begin{figure}[!ht]
	\begin{center}
		\includegraphics[width=6cm]{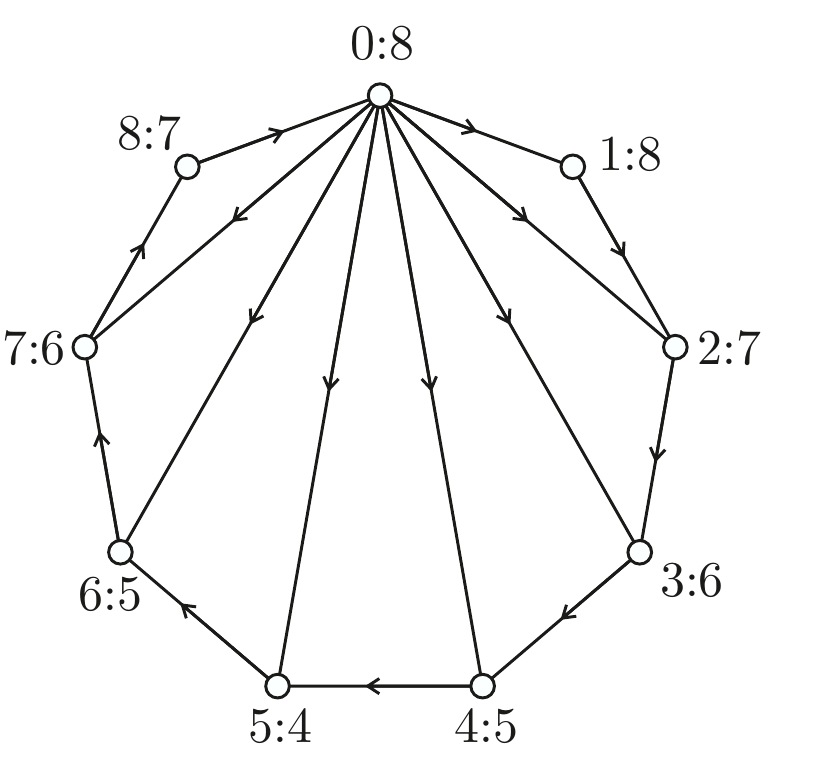}
	\end{center}
	\vskip-.5cm
	\caption{The strongly connected digraph  $C_9^*$ with inner out-radius 1 and inner in-radius 4. The label $i\!\!:\!\!j$ indicates that vertex $i$ has (inner) in-eccentricity $j$.}  
	\label{fig:Cn*}
\end{figure}
\end{proof}

We finish this section with a result on the mean inner distance.

\begin{Lemma}
	Let $X$ and $Y$ be random variables representing the inner distance between a pair of vertices randomly chosen in $G$ and in $LG$, respectively. Given a digraph $G$ with mean inner distance $\overline{d}_G=\mathbb{E}(X)$ (that is, the expectation of $X$), then the mean inner distance of $LG$ is
	$$
	\overline{d}_{LG}\leq \mathbb{E}(X)+1.
	$$
\end{Lemma}

\begin{proof}
	We compute $\overline{d}_{LG}$ as follows.
	\begin{align*}
	\overline{d}_{LG}&=\sum_{(u,v)\in U} \mathbb{E}(Y|u,v)\mathbb{P}(u,v) \leq \sum_{(u,v)\in U} (\dist(u,v)+1)\mathbb{P}(u,v) \\
	&=\sum_{(u,v)\in U} \dist(u,v)\mathbb{P}(u,v)+\sum_{(u,v)\in U} \mathbb{P}(u,v) = \mathbb{E}(X)+1,
	\end{align*}
	where  $U=\{(u,v):\dist(u,v)<\infty\}$, $\mathbb{E}(Y|u,v)$ is the conditional expectation of $Y$ assuming that the vertices $u$ and $y$ are taken, and $\mathbb{P}(u,v)$ is the probability of taking the vertices $u$ and $v$.
\end{proof}

\section{Sequences of inner diameters and orders}
\label{sec:sequences}

For a given digraph $G$, we are here interested in the sequences ${\cal D}=d_0,d_1,d_2,\ldots$ and ${\cal N}=n_0,n_1,n_2,\ldots$, where $d_k$ and $n_k$ are, respectively, the inner diameter and number of vertices  of the $k$-iterated line digraph $L^kG$, for $k=0,1,2,\ldots$,
where $L^0G=G$. 


Concerning ${\cal D}$, when the inner diameter coincides with the standard diameter, we already know the possible behaviors of $d_0,d_1,d_2,\ldots$, namely, $d_k$ is a constant if $G$ is a directed cycle or $d_k=d_{k-1}+1$ otherwise. Thus, we concentrate on the case when $G$ is not strongly connected.
In this case, the following result describes the different possibilities of the sequence of inner diameters $d_0=d(G), d_1=d(LG)$, $d_2=d(L^2G), \ldots$ depending on the existence of a strongly connected subdigraph $G'\subset G$. 

\begin{Proposition}
	\label{prop:diam-creix}
 Let $G$ be a non-strongly connected digraph. 
 \begin{itemize}
     \item[$(i)$]
If $G$ does not have any (non-trivial) strongly connected subdigraph $G'\subset G$, then there exists a value $h$ for which $d(L^{h}G)=0$ and  $L^{k}G=\emptyset$ when $k> h$.
\item[$(ii)$]
If $G$ has some (non-trivial) strongly connected subdigraph $G'\subset G$, there are two cases:
\begin{itemize}
   \item[$(ii.a)$] 
   If $G'$ is the unique such digraph and it is a directed cycle, then
   the sequence $d(G),d(LG),d(L^2G),\ldots$ becomes periodic from a certain iteration.
   \item[$(ii.b)$]
   If $G'$ is not a directed cycle, then the sequence $d(G),d(LG),d(L^2G),\ldots$ tends to infinity.
\end{itemize}
  \end{itemize}
\end{Proposition}
\begin{proof}
	$(i)$ Consider a digraph $G$ without any strongly connected subdigraph and its iterated line digraphs $LG, L^2G, L^3G, \ldots$ The sequence of their inner diameters can increase, decrease, or be constant at the beginning. Recall that a sequence is called unimodal if it first increases and then decreases. In any case, since no subdigraph of $G$ is strongly connected, then $G$ has a vertex $u$ with $|G^-(u)|=0$ (called a {\em source}) and a vertex $v$ with $|G^+(v)|=0$, (called a {\em sink}) where $u$ and $v$ are the first and the last vertices, respectively, of a longest path in $G$, say ${\cal P}:u_0(=u),u_1,\ldots,u_h(=v)$. Then, after in the iterated line digraph of $L^hG$, the vertex $u_0u_1\ldots u_h$ cannot be adjacent from or to any other vertex (otherwise, ${\cal P}$ would not be a longest path). Thus, since this applies to all paths of length $h$, $L^hG$ consists of some isolated vertices. Hence, $d(L^hG)=0$ and $L^kG=\emptyset$ for all $k>h$.\\
 $(ii.a)$ 
 Notice that if $G$ is periodic, this must be the case for the orders and the inner diameters of its iterated line digraphs. Then, the statement is a consequence of the results of Hemminger \cite{He74}, who, under the given hypothesis, determined the period of a finite periodic digraph.\\
 $(ii.b)$
 We already mentioned that if $G$ is a strongly connected digraph different from a directed cycle, then the inner diameter $d$ coincides with the standard diameter $D$ and, hence, $d(L^kG)=D(L^kG)=D(G)+k=d(G)+k$.
 \end{proof}
 


\begin{figure}[t]
	\begin{center}
		\includegraphics[width=14cm]{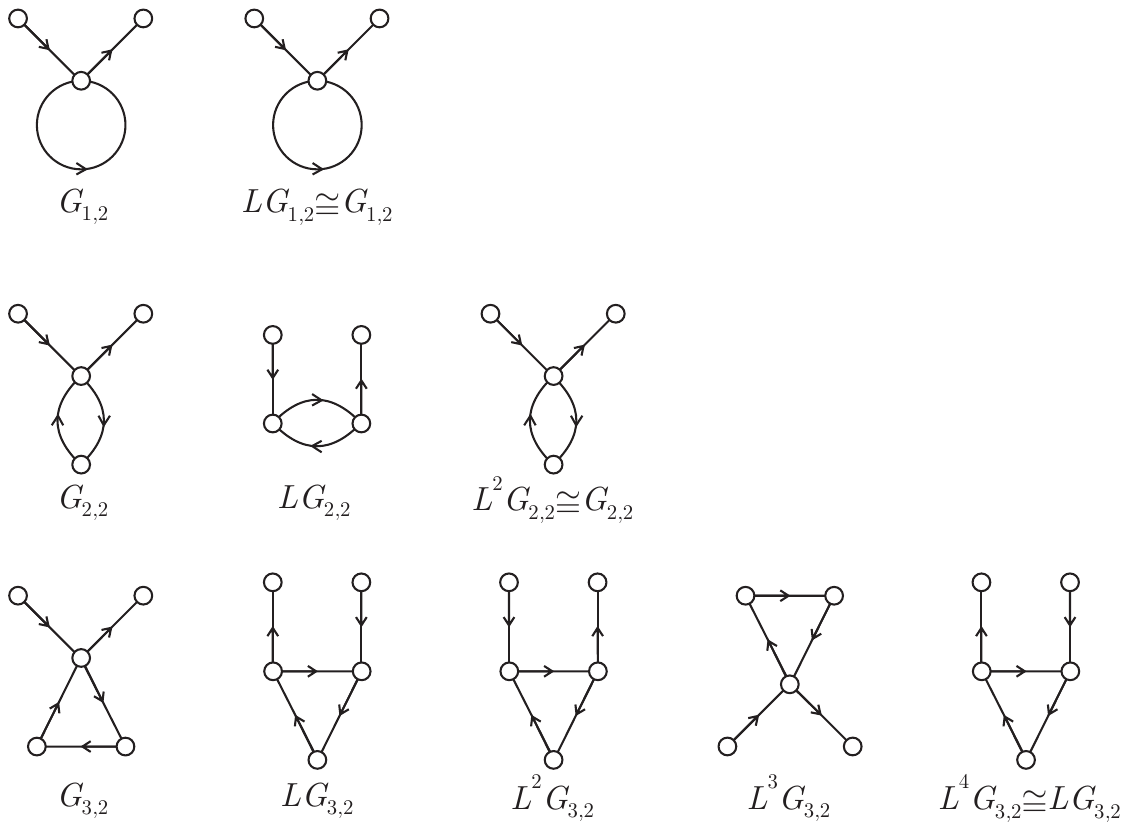}
	\end{center}
	\vskip-3.75cm
	\caption{$G_{1,2}$, $G_{2,2}$, and $G_{3,2}$ with some of its iterated line digraphs.}
	\label{fig:sequences}
\end{figure}

\begin{Example}
Here, we present a pair of examples regarding the cases of the last proposition.
\begin{itemize}
    \item[$(i)$] We can construct a digraph $G$ without any strongly connected subdigraph such that the sequence $d(G),d(LG),d(L^2G),\ldots$ increases by one unit at each iteration until any given value. After this value, the sequence decreases by one unit until zero. For example, we construct the K-shape digraphs denoted by $G_i$ for $i>2$. See the first six K-shaped digraphs in Figure \ref{fig:diam-creix}. In the sequence of the inner diameters of $G_i$ and its $\ell$-iterated line digraphs, the maximum inner diameter is $i+\ell_{\max}-1$, obtained for $\ell_{\max}=\lceil\frac{i}{2}-1\rceil$.
\begin{figure}[!ht]
	\begin{center}
		\includegraphics[width=12cm]{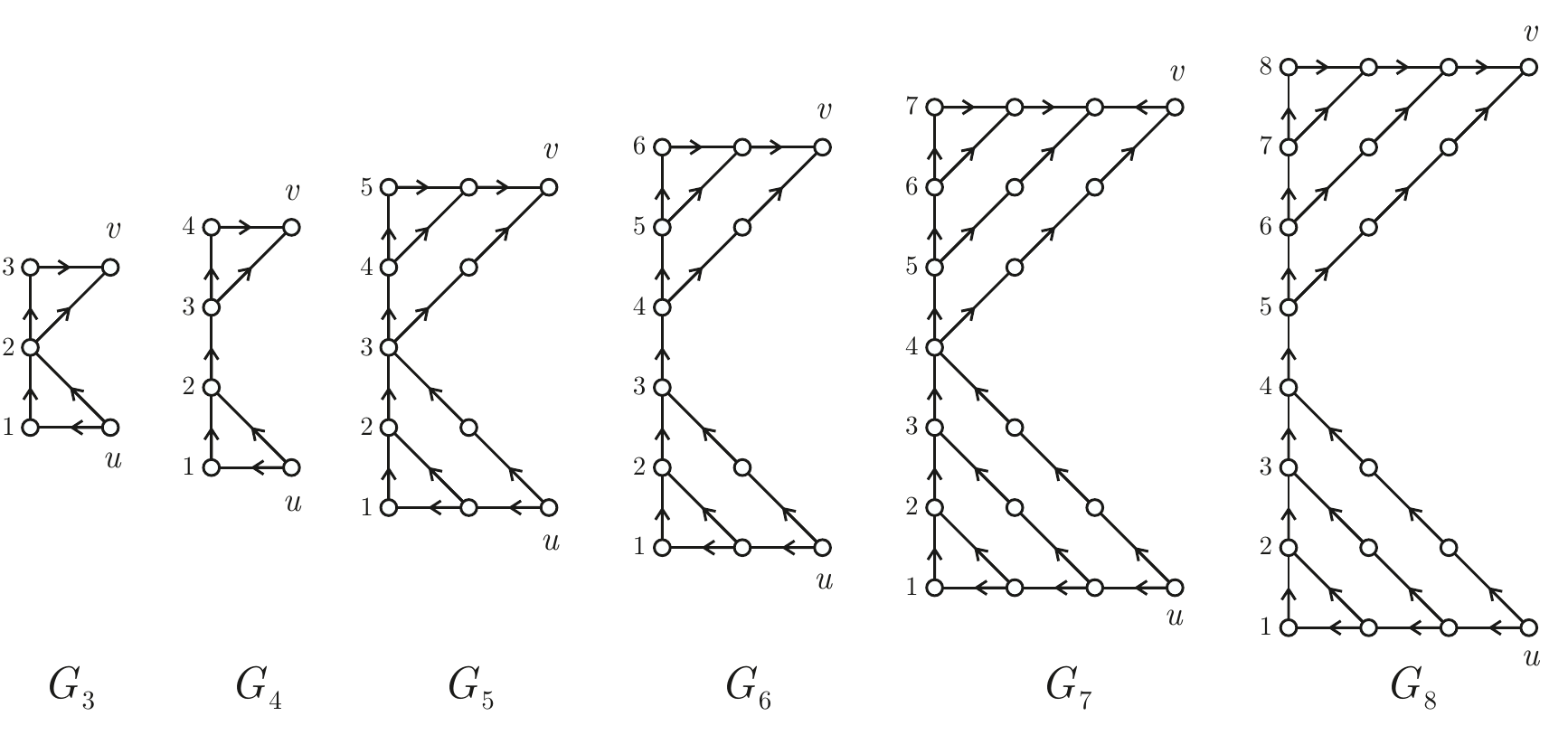}
	\end{center}
	\vskip-.5cm
	\caption{The first six K-shape digraphs, whose longest path is between vertices $u$ and $v$, with the following inner diameters:  $d(G_3)=2$, $d(G_4)=3$, $d(G_5)=4$, $d(G_6)=5$, $d(G_7)=6$, and $d(G_8)=7$.}
	\label{fig:diam-creix}
\end{figure}
    \item[$(ii.b)$] We construct the family of digraphs $G_{i,j}$ in Figure \ref{fig:sequences}, each of its digraphs formed by a directed cycle with a vertex with one out-adjacency and one in-adjacency.
\begin{itemize}
	\item[$\circ$] The inner diameters of $G_{1,2}$ and its iterated line digraphs are:\\
	 $2,2,2,2,\ldots$
	\item[$\circ$] The inner diameters of $G_{2,2}$ and its iterated line digraphs are:\\
	 $2,3,2,3,2,3,\ldots$
	\item[$\circ$] The inner diameters of $G_{3,2}$ and its iterated line digraphs are:\\
	 $3,4,3,3,4,3,3,4,\ldots$
	\item[$\circ$] The inner diameters of $G_{4,2}$ and its iterated line digraphs are:\\
	$4,5,4,4,4,5,4,4,4,5,\ldots$\\
	$\vdots$
	\item[$\circ$] The inner diameters of $G_{n,2}$ and its iterated line digraphs are:\\
	$n,n+1,n,\stackrel{(n-1)}{\ldots},$ $n,n+1,n,\stackrel{(n-1)}{\ldots},n,n+1,\ldots$
\end{itemize}
\end{itemize}
\end{Example}

Let us now consider the sequences ${\cal N}$.
In \cite{DaFi17}, Dalf\'o and Fiol gave a method to obtain the number of vertices of any iterated line digraph. 
In the following result, they obtained a recurrence equation on the number of
vertices $n_k$ of the
$k$-iterated line digraph of a digraph $G$.

\begin{Theorem}[\cite{DaFi17}]
	\label{maintheo}
    Let $G=(V,E)$ be a digraph on $n$ vertices, and consider a regular partition
	$\pi=(V_1,\ldots,V_m)$ with quotient matrix $\B$. Let $m(x)=x^r-\alpha_{r-1} x^{r-1}-\cdots-\alpha_0$ be the minimal polynomial of $\B$. Then, the number of vertices $n_k$ of the $k$-iterated line digraph $L^k(G)$ satisfies the
	recurrence
	\begin{equation*}
	\label{recur}
	n_k= \alpha_{r-1} n_{k-1}+\cdots+\alpha_{0} n_{k-r},\qquad \mbox{for } k=r,r+1,\ldots
	\end{equation*}
	initialized with the values $n_{k}$, for $k=0,1,\ldots,r-1$, given by 
	\begin{equation}
	\label{n_l}
	n_{k}=\sum_{i=1}^m |V_i|\sum_{j=1}^m (\B^{k})_{ij}=\s\B^{k}\j^{\top},
	\end{equation}
	where $\s=(|V_1|,\ldots,|V_m|)$  and $\j=(1,\ldots,1)$.
\end{Theorem}

This result allows us to give a method to find the number of words of length $n$ over an alphabet of $d$
symbols, avoiding a given set $S$ of subwords.
Basically, the method consists of the following steps:
\begin{enumerate}
    \item 
    Take the De Bruijn digraph $B(d,n')$, where $n'$ is the maximum length of the forbidden subwords in $S$.
    \item 
    Obtain the digraph $G$ obtained from $B(d,n')$ by removing all vertices containing some subwords in $S$.
    \item 
    Compute the minimal polynomial of the adjacency matrix $\A$ of $G$.
    \item 
    Apply Theorem~\ref{maintheo} to obtain the recurrence formula to obtain the numbers of vertices of the iterated line digraphs of $G$. 
 \end{enumerate}

\label{sec:examples}
To illustrate our method, 
in what follows, we give examples of the four possible behaviors of the sequence
$n_0, n_1, n_2,\ldots$ of the number of vertices (or inner diameter) of the iterated line digraphs of a given digraph.
Namely, when it is increasing, periodic, tending to a positive constant, or tending to zero.

\subsection{Binary sequences}

Suppose we want to know the number $n_k$ of binary words of length $k$ that do not contain the subwords $000$, $001$, $010$, and $100$.
For $k\ge 3$, this is the number of vertices of the iterated line digraph $L^{k-3}G$, where $G$ is the digraph obtained from the De Bruijn digraph $B(2,3)$ by deleting the vertices $000$, $001$, $010$, and $100$, as shown in Figure \ref{fig:B'(2,3)}.
\begin{figure}[!ht]
	\begin{center}
		\includegraphics[width=6cm]{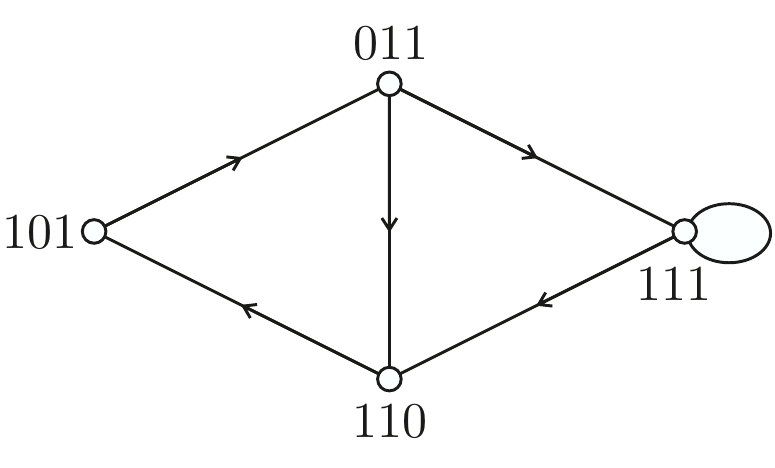}
	\end{center}
	\vskip-.5cm
	\caption{$B(2,3)$ minus $000$, $001$, $010$, and $100$.}
	\label{fig:B'(2,3)}
\end{figure}

Then, the adjacency matrix of $G$ is
$$
\A=\left(
\begin{array}{cccc}
1& 1& 0& 0\\
0& 0& 1& 0\\
0& 0& 0& 1\\
1& 1& 0& 0    
\end{array}\right),
$$
with minimal polynomial $m(x)=x^4-x^3-x$. Thus, by Theorem~\ref{maintheo}, the number of vertices of $L^k(G)$ satisfies the recurrence $n_k=n_{k-1}+n_{k-3}$ for $k\ge 3$ and we obtain the sequence
$$
4, 6, 9, 13, 19, 28, 41, 60, 88, 129, 189, 277, 406, 595, 872, 1278,\ldots
$$
From $n\ge 5$, this corresponds to the sequence $a(n)$ for $n=0,1,2,\ldots$ referred as $A000930$ in OEIS \cite{Sl}, and named {\em Narayana's cows sequence}, with initial terms $a(0)=a(1)=a(2)=1$ and recurrence $a(n) = a(n-1) + a(n-3)$.
Among other possible interpretations, and in concordance with our results, in \cite{Sl}, it is commented that $a(n+2)$ is the number of $n$-bit $0$-$1$ sequences that avoid both $00$ and $010$ (Callan, 2004).

In Table \ref{table:forbid-words}, we show other possibilities when we remove three vertices of $B(2,3)$.

In the Appendix, we show tables with examples of new sequences (similar to Table \ref{table:forbid-words}), which are not in OEIS. The new sequences were obtained starting from De Bruijn and Square-free digraphs.

\subsubsection{Cyclic Kautz digraphs}
\label{sec:ex-CK}

The {\em cyclic Kautz digraph} $CK(d,\ell)$, introduced by B\"{o}hmov\'{a},
Dalf\'{o}, and Huemer in~\cite{BoDaHu17}, has vertices labeled by all possible sequences $a_1a_2\ldots a_\ell$  with $a_i\in\{0,1,\ldots,d\}$, $a_i\neq a_{i+1}$ for $i=1,\ldots,\ell-1$, and $a_1\neq a_\ell$. Moreover,
there is an arc from vertex $a_1 a_2\ldots a_\ell$ to vertex $a_2 \ldots a_\ell
a_{\ell+1}$, whenever $a_{\ell+1}\neq a_{\ell},a_2$. 

\begin{table}[!ht]
{\small
\begin{center}
\begin{tabular}{|c|l|l|}
\hline
\mbox{Forbidden subwords} & Sequence  & OEIS \cite{Sl}   \\
\hline
000, 010, 011 & 5, 5, 5, 5, 5, 5, 5, 5, \ldots
 & A010716\tablefootnote{Constant sequence: the all 5's sequence. } for $n\ge 0$ \\
000, 001, 101 &  5, 6, 6, 6, 6, 6, 6, 6, \ldots &
A101101\tablefootnote{$a(1)=1$, $a(2)=5$, and $a(n)=6$ for $n\ge3$} for $n\geq 2$ \\
001, 010, 011 & 5, 6, 7, 8, 9, 10, 11, 12, \ldots &
A000027\tablefootnote{The positive integers.} 	for $n\geq 5$ \\
000, 010, 111 & 5, 6, 7, 9, 11, 13, 16, 20, \ldots
&
A164317\tablefootnote{Number of binary strings of length $n$ with no substrings equal to 000, 010, or 111.} for $n\geq 3$ \\
000, 011, 110 & 
5, 6, 8, 10, 13, 17, 22, 29, \ldots
& A052954\tablefootnote{Expansion of $\frac{2-x-x^2-x^3}{(1-x)(1-x^2-x^3)}$.} for $n\geq 8$	\\
000, 010, 101 & 5, 7, 10, 14, 19, 26, 36, 50, \ldots
& A003269\tablefootnote{ $a(n) = a(n-1) + a(n-4)$ with $a(0) = 0$, $a(1) = a(2) = a(3) = 1$.}  for $n\geq 8$	\\
001, 010, 100 & 5, 7, 10, 14, 20, 29, 42, 61, \ldots
& A020711\tablefootnote{Pisot sequences $E(5,7)$, $P(5,7)$.}  for $n\geq 0$ \\
000, 001, 010	& 5, 7, 11, 16, 23, 34, 50, 73, \ldots
& A164316\tablefootnote{Number of binary strings of length $n$ with no substrings equal to 000, 001, or 010. }	for $n\geq 3$\\
000, 001, 011 & 5, 7, 8, 10, 11, 13, 14, 16, \ldots
& A001651\tablefootnote{Numbers that are not divisible by 3.}  for $n\geq 3$ \\
001, 010, 101 &	5, 7, 9, 11, 13, 15, 17, 19, \ldots
& A005408\tablefootnote{The odd numbers: $a(n) = 2n + 1$.}  for $n\geq 2$. \\
000, 001, 111 & 5, 7, 9, 12, 16, 21, 28, 37, \ldots 
& A000931\tablefootnote{Padovan sequence (or Padovan numbers): $a(n) = a(n-2) + a(n-3)$ with $a(0) = 1$, $a(1) = a(2) = 0$.}
 for $n\geq 12$ 
 \\
000, 001, 100 & 5, 8, 13, 21, 34, 55, 89, 144, \ldots
& A000045\tablefootnote{Fibonacci numbers: $F(n) = F(n-1)+F(n-2)$ with $F(0) = 0$ and $F(1) = 1$.} for $n\geq 5$	\\
\hline
\end{tabular}
\caption{Forbidden words in the De Bruijn digraphs and the sequence obtained with the numbers of vertices on $L^0(G)=G, L^1(G), L^2(G),\ldots$}
\label{table:forbid-words}
\end{center}
}
\end{table}

\begin{figure}[t]
    \vskip-.5cm
    \begin{center}
        \includegraphics[width=14cm]{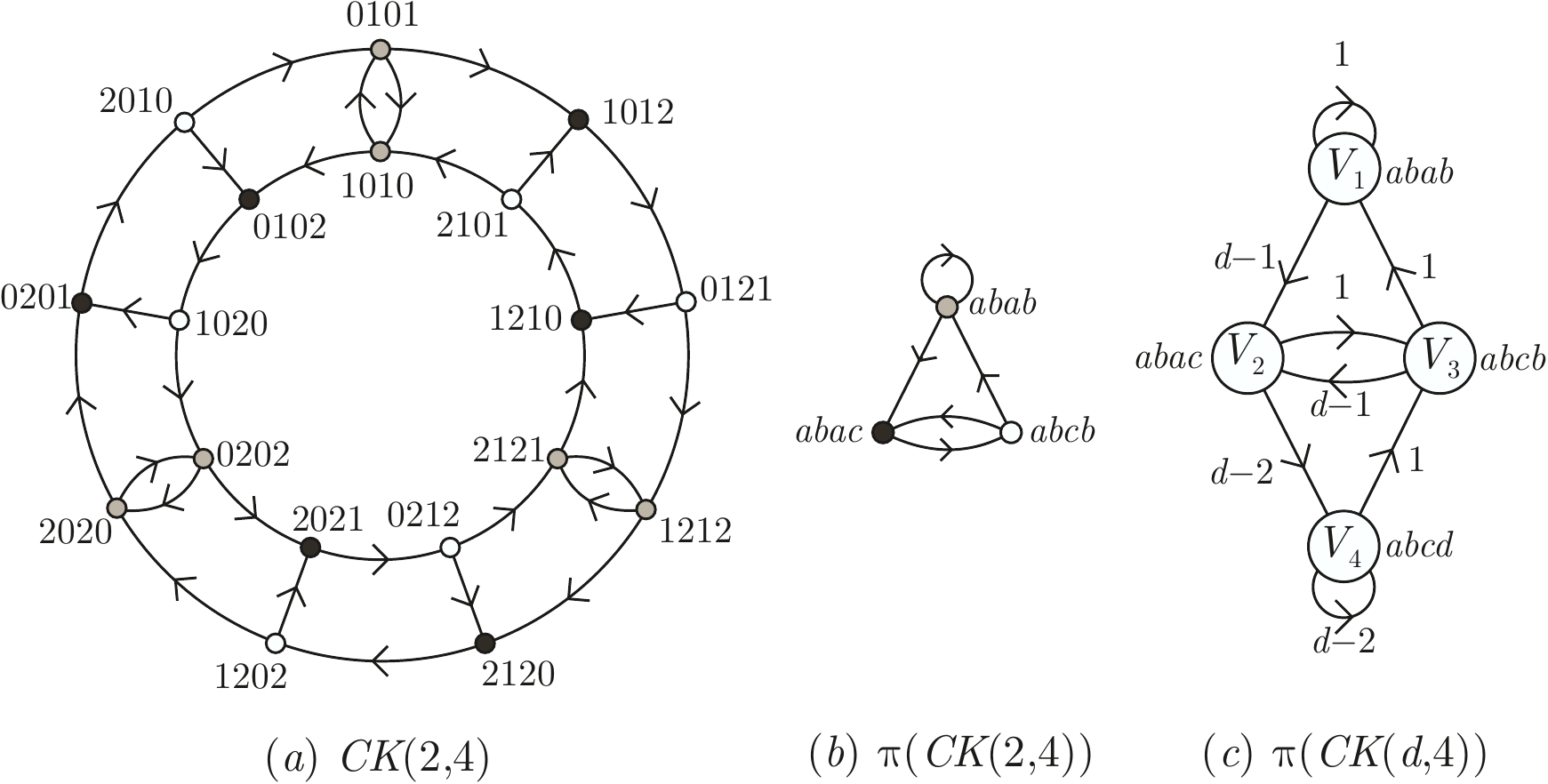}
    \end{center}
	\caption{$(a)$ The cyclic Kautz digraph $CK(2,4)$; $(b)$ its quotient digraph $\pi(CK(2,4))$; and
$(c)$ the quotient digraph of $CK(d,4)$.}
	\label{fig:quocient-donut-color}
\end{figure}


\begin{figure}[t]
    \vskip-.5cm
    \begin{center}
        \includegraphics[width=6cm]{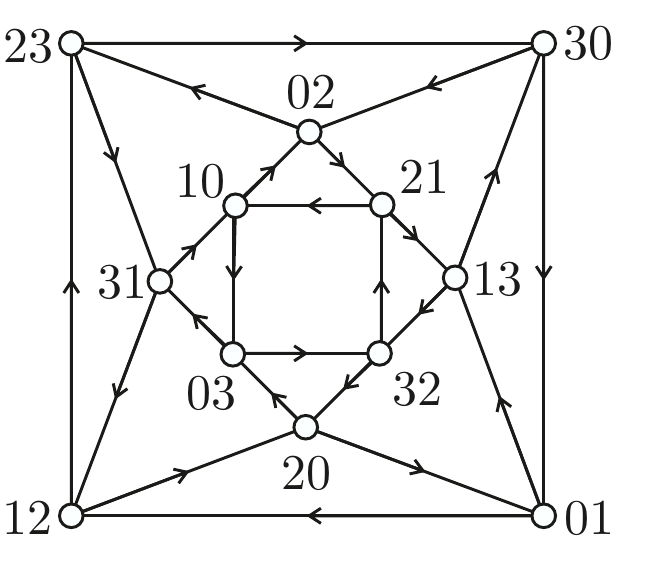}
         \includegraphics[width=7cm]{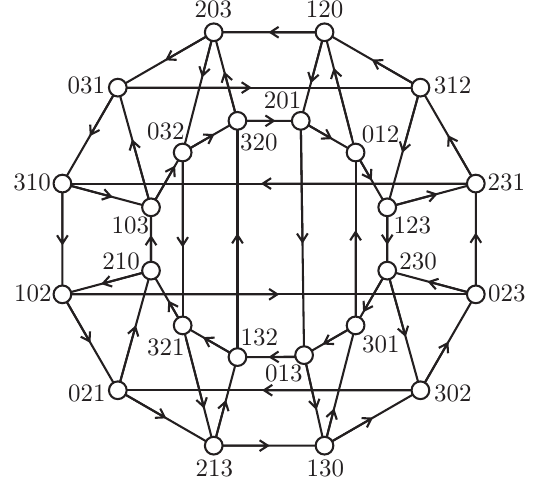}
    \end{center}
	\caption{The subKautz digraph $sK(3,2)$ (left) and its line digraph $CK(3,3)=L(sK(3,2))$ (right).}
	\label{fig:K'(3,2)&CK(3,3)}
\end{figure}

For example, Figure~\ref{fig:quocient-donut-color}$(a)$ shows the cyclic Kautz digraph $CK(2,4)$.
By this definition, we observe that the cyclic Kautz digraph $CK(d,\ell)$ is a subdigraph of the well-known Kautz digraph $K(d,\ell)$, that has vertices $a_1a_2\ldots a_\ell$ with $a_i\neq a_{i+1}$ for $i=1,\ldots,\ell-1$, (so,  without the requirement $a_1\neq a_\ell$). Thus, there is an arc from vertex $a_1 a_2\ldots a_\ell$ to vertex $a_2 \ldots a_\ell
a_{\ell+1}$, whenever $a_{\ell+1}\neq a_{\ell}$.

Moreover, the cyclic Kautz digraph is the line digraph of the so-called subKautz digraph (see Dalf\'o \cite{Da17(b)}) defined as follows.
Given integers $d$ and $\ell$, with $d,\ell\ge 2$, a {\em subKautz digraph} $sK(d, \ell)$ has the same set of vertices as $K(d,\ell)$, and adjacencies
$a_1a_2\ldots a_\ell \rightarrow a_2\ldots a_\ell a_{\ell+1}$, with $a_{\ell+1} \neq a_1, a_{\ell}$.
Then,  the out-degree of a vertex $a_1a_2\ldots a_\ell$ is $d$ if $a_1 = a_{\ell}$, and $d-1$ otherwise. In
particular, the subKautz digraph $sK(d, 2)$ is $(d-1)$-regular and it can be obtained from the
Kautz digraph $K(d, 2)$ by removing all its arcs forming a digon.

\begin{Lemma}[\cite{Da17(b)}]
 The cyclic Kautz digraph $CK(d,\ell)$ is isomorphic to the line digraph of the subKautz digraph $sK(d,\ell-1)$, that is,   $CK(d,\ell)=$ \linebreak $L(sK(d,\ell-1))$. 
\end{Lemma}

For example, the subKautz digraph $sK(3,2)$ and the cyclic Kautz digraph $CK(3,3)=L(sK(3,2))$ are shown in Figure
\ref{fig:K'(3,2)&CK(3,3)}.

Since, in general,  the subKautz and cyclic Kautz digraphs are not $d$-regular, the number of vertices of their iterated line digraphs are not obtained by repeatedly multiplying by $d$. Instead, we can apply our method, as shown in what follows with the cyclic Kautz digraph $CK(2,4)$. This digraph has a regular partition $\pi$ of its vertex set into three classes (each one with 6 vertices): $abcb$ (the second and the last digits are equal), $abab$ (the first and the third digits are equal, and also the second and the last), and $abac$ (the first and the third digits are equal). Then, the quotient matrix of $\pi$, which coincides with the adjacency matrix of $\pi(CK(2,4))$, is 
$$
\B=\left(
  \begin{array}{ccc}
    0 & 1 & 1 \\
    0 & 1 & 1 \\
    1 & 0 & 0 \\
  \end{array}
\right),
$$
where the order of the vertices is $abcb$, $abab$, $abac$ 
and it has minimal polynomial $m(x)=x^3-x^2-x$. Consequently, by
Theorem~\ref{maintheo}, the number of vertices of $L^k(CK(2,4))$ satisfies the recurrence $n_k=n_{k-1}+n_{k-2}$ for $k\ge 3$. In fact, in this
case, $\s(\B^2-\B-\I)\j^{\top}=0$, and the above recurrence applies from $k=2$.
This, together with the initial values $n_0=18$ and $n_1=\s\B\j^{\top}=30$, yields the Fibonacci sequence, $n_2=48, n_3=78, n_4=126, n_5=204,\ldots$,
as B\"{o}hmov\'{a}, Dalf\'{o}, and Huemer~\cite{BoDaHu17} proved by using a
combinatorial approach.

Curiously enough, $n_k$ is also the number of ternary length-2 square-free words of length
$k+4$; see the sequence A022089 in OEIS (The On-Line Encyclopedia of Integer Sequences) \cite{Sl}. To prove this, we suggest working with the square-free digraphs, as shown in the following subsection. 

In fact, our method allows us to generalize this result and, for instance, derive a
formula for the order of $L^k(CK(d,4))$ for any value of the degree $d\ge 2$. To
this end, it is easy to see that a quotient digraph of $CK(d,4)$ for $d>2$ is as
shown in Figure~\ref{fig:quocient-donut-color}$(c)$, where now we have to
distinguish four classes of vertices.
Then, the corresponding quotient matrix is
$$
\B=\left(
  \begin{array}{cccc}
    1 & d-1 & 0 & 0 \\
    0 & 0   & 1 & d-2\\
    1 & d-1 & 0 & 0 \\
    0 & 0   & 1 & d-2
  \end{array}
\right),
$$
and it has minimal polynomial is $m(x)=x^3-(d-1)x^2-x$.
In turn, this leads to the recurrence formula $n_k=(d-1)n_{k-1}+n_{k-2}$, with
initial values $n_0=d^4+d$ and $n_1=d^5-d^4+d^3+2d^2-d$, which are
computed by using \eqref{n_l} with the vector
\begin{align*}
  \s & =(|V_1|,|V_2|,|V_3|,|V_4|) \\
  & =((d+1)d, (d+1)d(d-1), (d+1)d(d-1), (d+1)d(d-1)(d-2)).
\end{align*}
Solving the recurrence, we get the closed formula
$$
n_k =\frac{2^k d}{\sqrt{\Delta}}
\left(\frac{
(d^2+d)\sqrt{\Delta}-d^3-d-2
}
{(1-d-\sqrt{\Delta})^{k+1}}
+\frac{(d^2+d)\sqrt{\Delta}+d^3+d+2
}
{(1-d+\sqrt{\Delta})^{k+1}}\right),
$$
where $\Delta=d^2-2d+5$ and, hence, $n_k$ is an increasing sequence. For instance, with $d=3$ (four symbols), we get the sequence $84,204,492,1188,2868,\ldots$ Dividing the terms by 12 we get $7,17,41,99,239,\ldots$ which correspond to the sequence $A001333(n)$ for $n=3,4,\ldots$ in OEIS~\cite{Sl}.

\subsubsection{Square-free digraphs}
\label{sec:ex-SF}

\begin{figure}[t]
    \vskip-.5cm
    \begin{center}
        \includegraphics[width=14cm]{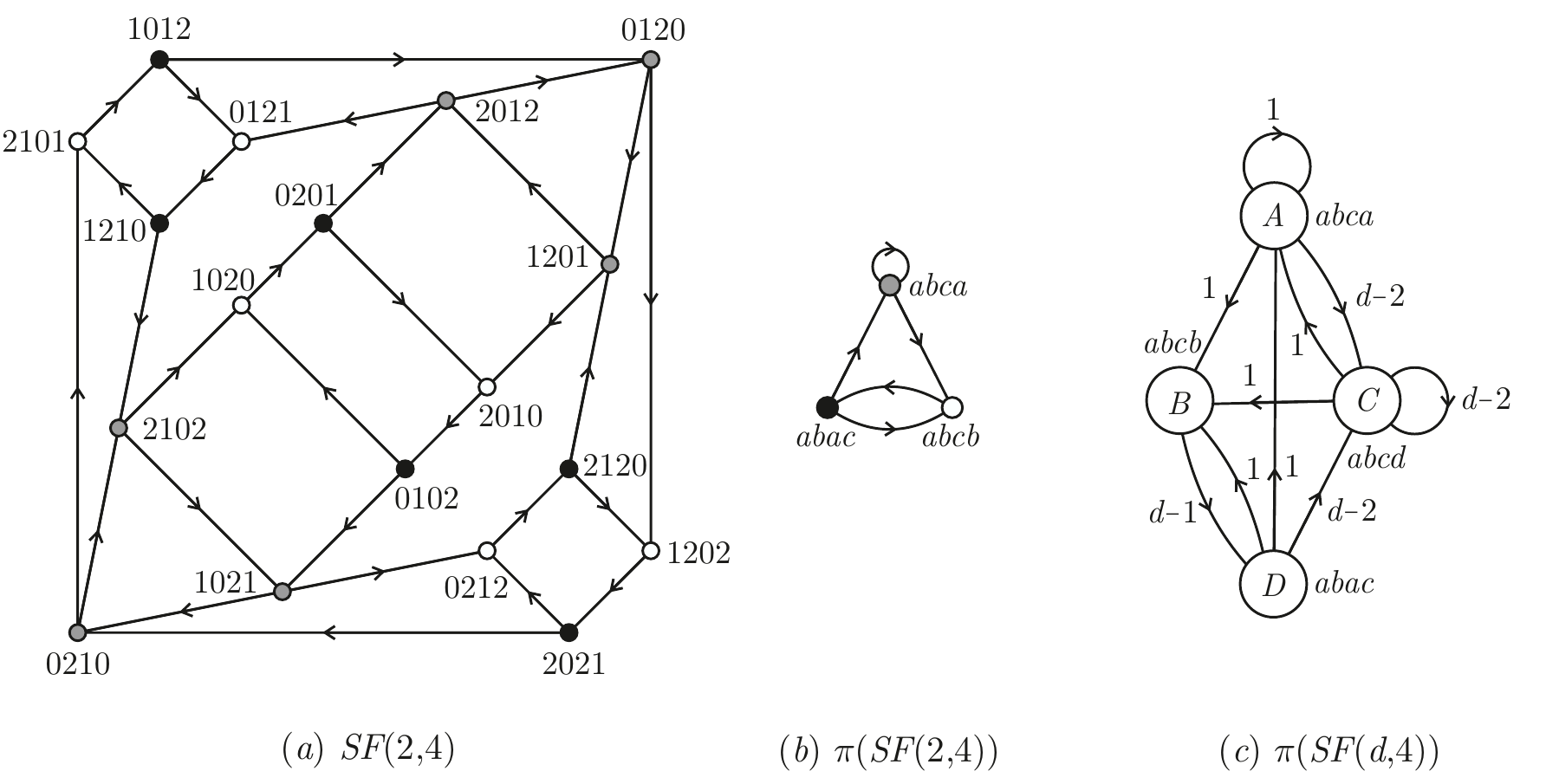}
    \end{center}
    \vskip-.5cm
	\caption{$(a)$ The square-free digraph $SF(2,4)$; $(b)$ its quotient digraph $\pi(SF(2,4))$; and
$(c)$ the quotient digraph of $SF(d,4)$.}
	\label{fig:SF}
\end{figure}

The {\em square free digraph} $SF(d,\ell)$ has vertices labeled with words $a_1a_2\ldots a_{\ell}$, on a $(d+1)$-letter alphabet that does not contain an adjacent
repetition of any subword of length at least $2$. Moreover, there is an arc from $a_1a_2\ldots a_{\ell}$ to $a_2\ldots a_{\ell}a_{\ell+1}$ when $a_{\ell+1}\neq a_{\ell}$ and, if $a_{\ell-2}=a_{\ell}$, then $a_{\ell+1}\neq a_{\ell-1}$ (to avoid equal consecutive subwords of length two). 
Notice that $SF(d,\ell)$ can be obtained from the De Bruijn digraph $B(d+1,\ell)$ by removing the vertices with the forbidden labels. An example is the square free digraph $SF(2,4)$ shown in Figure \ref{fig:SF}$(a)$.
Now, the reason for $L^k(CK(2,4))$ and $L^k(SF(2,4))$ sharing the same number of vertices for every $k\ge 0$ is that they have the same quotient digraphs shown in Figures  \ref{fig:quocient-donut-color}$(b)$ and \ref{fig:SF}$(b)$, respectively.
Let us now derive a
formula for the order of $L^k(SF(d,4))$ for any value of the degree $d\ge 2$. To
this end, it is easy to see that a quotient digraph of $SF(d,4)$ for $d>2$ is as
shown in Figure~\ref{fig:SF}$(c)$, where now we have to
distinguish four classes of vertices.
Then, the corresponding quotient matrix is
$$
\B=\left(
  \begin{array}{cccc}
    1 & 1 & d-2 & 0\\
    0 & 0 & 0   & d-1\\
    1 & 1 & d-2 & 0\\
    1 & 1 & d-2 & 0 
  \end{array}
\right),
$$
and it has minimal polynomial $m(x)=x^3-(d-1)x^2-(d-1)x$.
In turn, this leads to the recurrence formula $n_k=(d-1)n_{k-1}+(d-1)n_{k-2}$, with
initial values $n_0=d^4+d^3-d^2-d$ and $n_1=d^5+d^4-2d^3-d^2+d$, which are
computed by using \eqref{n_l} with the vector
\begin{align*}
  \s & =(|V_1|,|V_2|,|V_3|,|V_4|) \\
     & =((d^2 - 1)d, (d^2 - 1)d, (d^2 - 1)d(d - 2), (d^2 - 1)d).
\end{align*}
For $d=2$, we obtain, as expected, the same sequence $18,30,48,78,\ldots$ for the number of vertices of $L^k(SF(2,4))$ as in $L^k(CK(2,4))$. But this does not occur for other values of $d$. For instance, for $d=3$, the numbers of $L^k(SF(3,4))$ follow the sequence $96,264,720,1968,5376,14688,\ldots$, which correspond to the values $a(n)$, for $n\ge 4$, of A239171 in  \cite{Sl}. This is a special case, for $k=1$, of A239178. There, $T(n,k)$ is the number of $(n+1)\times (k+1)$ $(0,1,2)$-arrays with no element greater than all horizontal neighbors or equal to all vertical neighbors. In our case, $a(n)=T(n,1)$ can also be defined as the number of ternary words of length $n+1$ with no three consecutive equal symbols and no first or last two equal symbols. For instance, for $n=4$, the 16 words beginning with $01$ are
\begin{equation*}
\begin{tabular}{cccccccc}
01001 & 01002 & 01010 & 01012 & 01020 & 01021 & 01101 & 01102\, \\
01120 & 01121 & 01201 & 01202 & 01210 & 01212 & 01220 & 01221.
\end{tabular}
\end{equation*}
Thus, since we have the same number for all 6 possible first two (different) symbols, we have a total of $a(4)=6\cdot 16=96$, as shown above.

\begin{figure}[t]
    \vskip-.5cm
    \begin{center}
        \includegraphics[width=14cm]{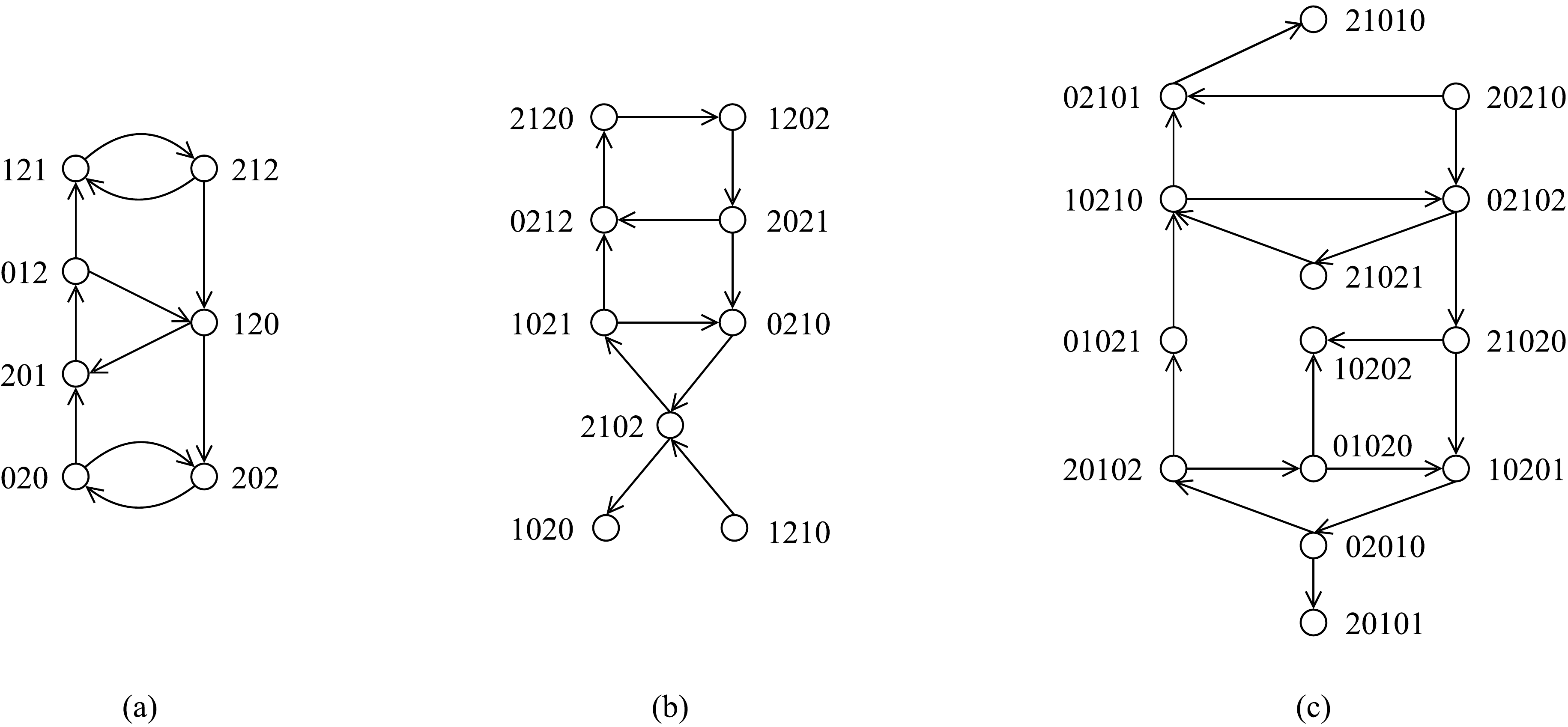}
    \end{center}
    \vskip-.5cm
	\caption{(a) The square-free digraph $SF'(2,3)$ with forbidden subsequences $021$ and $10$ in its vertex labels. (b) The square-free digraph $SF'(2,4)$ with forbidden subsequence $01$. (c) The square-free digraph $SF'(2,5)$ with forbidden subsequence $12$.}
	\label{fig:SF_example}
\end{figure}


We found some unique sequences for the number of vertices in $L^k(SF(d,l))$ that are not listed in \cite{Sl}. In Figure \ref{fig:SF_example}, we present three square-free digraphs with forbidden sub-sequences in vertex labels that produce these unique sequences. The first square-free digraph $SF(2,3)$ (Figure \ref{fig:SF_example}(a)) was constructed by forbidding sub-sequences $021$ and $10$ in vertex labels and produces the unique sequence $7, 11, 16, 24, 36, 53, 80, 118, 177, 263,\ldots$ for the number of vertices in $L^k(SF(2,3))$. We forbid the subsequence $01$ in the vertex labels of the square-free digraph $SF(2,4)$ (Figure \ref{fig:SF_example}(b)). The unique sequence for the number of vertices in $L^k(SF(2,4))$ is $9, 11, 14, 17, 20, 25, 30, 37, 45, 55,\ldots$ The square-free digraph $SF(2,5)$ (Figure \ref{fig:SF_example} (c)) cannot have the subsequence $12$ in the vertex labels and generates the unique sequence $14, 18, 22, 27, 32, 40, 48$, $59, 72, 88,\ldots$ for the number of vertices in $L^k(SF(2,5))$.

\subsubsection{Unicyclic digraphs}
\label{sec:ex-unicyclic}

\begin{figure}[t]
    \begin{center}
        \includegraphics[width=10cm]{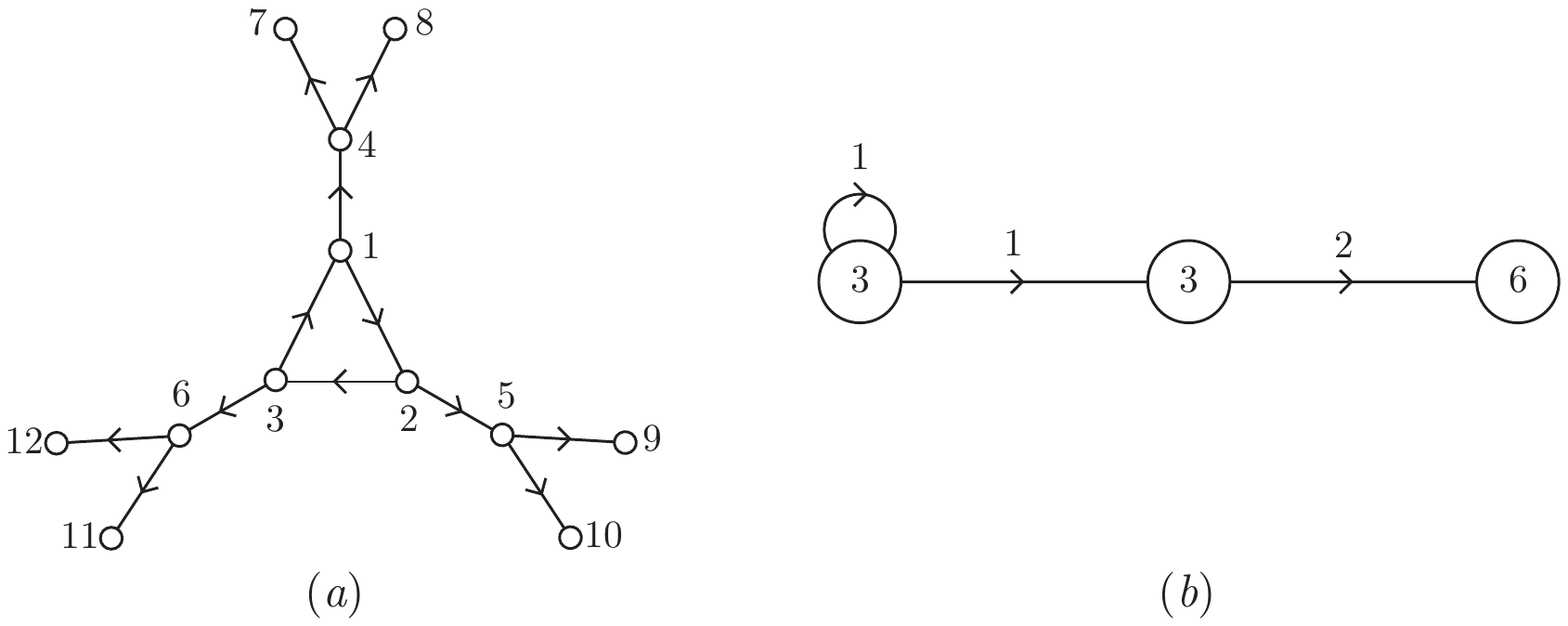}
    \end{center}
    \vskip-3.5cm
	\caption{The unicyclic digraph $G_{3,2}$ and its quotient digraph.}
	\label{fig:unicyclic}
\end{figure}

A unicyclic digraph is a digraph with exactly one (directed) cycle. As usual, we
denote a cycle on $n$ vertices by $C_n$. For example, consider the digraph
$G_{n,d}$, obtained by joining to every vertex of $C_n$ one `out-tree' with $d$
leaves (or `sinks'), as shown in Figure~\ref{fig:unicyclic}$(a)$ for the case
$G_{3,2}$. This digraph has the regular partition $\pi=(V_1, V_2, V_3)$, where $V_1$ is the set of vertices of the cycle, $V_2$ the central vertices of the trees,  and $V_3$ the set of leaves. (In the figure $V_1=\{1,2,3\}$, $V_2=\{4,5,6\}$, and $V_3=\{7,8,9,10,11,12\}$).
This partition gives the quotient digraph $\pi(G)$ of
Figure~\ref{fig:unicyclic}$(b)$, and the
quotient matrix
$$
\B=\left(
  \begin{array}{ccc}
    1 & 1 & 0 \\
    0 & 0 & d \\
    0 & 0 & 0 \\
  \end{array}
\right),
$$
with minimal polynomial $m(x)=x^3-x^2$. Then, by Theorem \ref{maintheo}, the order
of $L^k(G)$ satisfies the recurrence $n_k=n_{k-1}$ for $k\ge 1$, since
$\s(\B^{k}-\B^{k-1})\j^{\top}$ $=0$ for $k=1,2$, where $\s=(n,n,nd)$. Thus, we conclude that all the iterated line digraphs $L^k(G)$ have constant order
$n_k=n_0=n(d+2)$, that is, $n_k$ tends to a positive constant.
(In fact, this is because in this case $L(G)$---and, hence, $L^k(G)$---is
isomorphic to $G$.)

\subsubsection{Acyclic digraphs}
\label{sec:ex-acyclic}

\begin{figure}[t]
    \begin{center}
        \includegraphics[width=9cm]{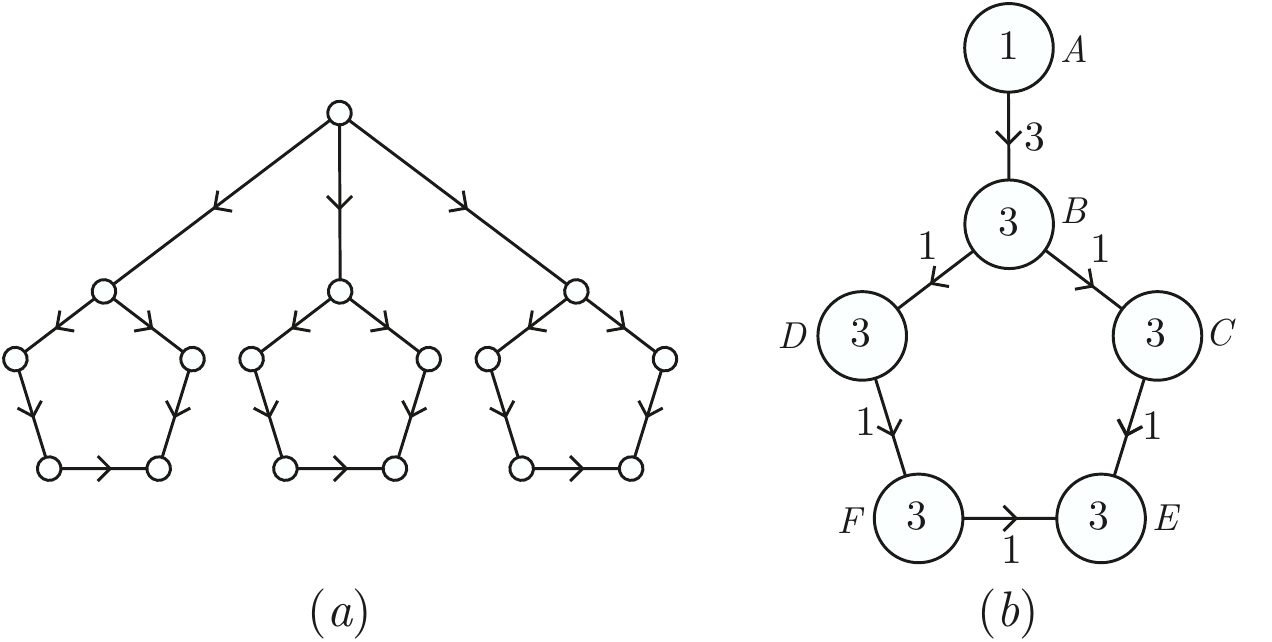}
    \end{center}
  \caption{An acyclic digraph and its quotient digraph.}
	\label{fig:acyclic}
\end{figure}

Finally, let us consider an example of an acyclic digraph, that is, a digraph
without directed cycles, such as the digraph $G$ of Figure~\ref{fig:acyclic}$(a)$.
Its quotient digraph is depicted in Figure~\ref{fig:acyclic}$(b)$, with quotient
matrix
$$
\B=\left(
  \begin{array}{cccccc}
    0 & 3 & 0 & 0 & 0 & 0 \\
    0 & 0 & 1 & 1 & 0 & 0 \\
    0 & 0 & 0 & 0 & 1 & 0 \\
    0 & 0 & 0 & 0 & 0 & 1 \\
    0 & 0 & 0 & 0 & 0 & 0 \\
    0 & 0 & 0 & 0 & 1 & 0 \\
  \end{array}
\right),
$$
and minimal polynomial $m(x)=x^{5}$. This indicates that $n_k=0$ for every $k\ge 5$
(as expected, because $G$ has no walk of length larger than or equal to $5$). Moreover, from \eqref{n_l}, the first values are
$n_0=16$, $n_1=18$, $n_2=15$, $n_3=9$, and $n_4=3$.


\subsection*{Acknowledgment}
\thanks{
This research has been supported by AGAUR from the Catalan Government under project 2021SGR00434 and MICINN from the Spanish Government under project PID2020-115442RB-I00.
The research of M. A. Fiol was also supported by a grant from the  Universitat Polit\`ecnica de Catalunya with references AGRUPS-2022 and AGRUPS-2023. The research of D. Z\'avack\'a was supported by G-24-158-00 and VEGA 1/0437/23.}

\section*{Declaration of competing interest}
The authors have no competing interests.



\newpage
\section*{Appendix}
Here, we show tables with examples of new sequences (similar to Table \ref{table:forbid-words}) that are not in OEIS. The new sequences were obtained starting from De Bruijn and Square-free digraphs.

\begin{table}[!ht]
\scriptsize
\begin{center}
\begin{tabular}{|c|l|l|}
\hline
Forbidden subwords & Sequence & OEIS  \\
\hline
000, 1001, 1011, 1101 & 10, 13, 13, 14, 13, 14, 13, 14, 13, 14, 13, 14, 13, \ldots & not in OEIS\\
0001, 0010, 0110, 111 & 10, 13, 14, 16, 17, 19, 20, 22, 23, 25, 26, 28, 29, \ldots & not in OEIS\\
0010, 0101, 0110, 111 & 10, 13, 15, 17, 18, 18, 18, 18, 18, 18, 18, 18, 18,  \ldots & not in OEIS\\

0000, 0101, 0110, 111 & 10, 13, 15, 19, 23, 28, 34, 42, 51, 62, 76, 93, 113,  \ldots & not in OEIS\\
000, 0111, 1010, 1011 & 10, 13, 16, 21, 27, 33, 41, 52, 64, 78, 97, 120, 146, \ldots & not in OEIS\\
000, 0111, 1100, 1101 & 10, 13, 16, 21, 27, 35, 46, 60, 79, 104, 137, 181, \ldots & not in OEIS\\
0000, 0100, 0101, 111 & 10, 13, 16, 22, 30, 39, 51, 68, 91, 120, 158, 210, \ldots & not in OEIS\\
0010, 0101, 110, 1101 & 10, 13, 17, 21, 25, 29, 33, 37, 41, 45, 49, 53,  \ldots & not in OEIS\\
000, 0110, 1001, 1101 & 10, 13, 17, 22, 26, 31, 35, 40, 44, 49, 53, 58, \ldots & not in OEIS\\
000, 0101, 110, 1101 & 10, 13, 17, 22, 28, 35, 44, 55, 68, 84, 104, 128,  \ldots & not in OEIS\\
000, 0100, 0111, 1101 & 10, 13, 17, 22, 29, 37, 48, 61, 79, 100, 129, 163,  \ldots & not in OEIS\\
000, 0101, 0111, 1100 & 10, 13, 17, 23, 29, 37, 49, 61, 77, 101, 125, 157, \ldots & not in OEIS\\
000, 0111, 1010, 1101 & 10, 13, 17, 23, 31, 41, 54, 71, 93, 122, 160, 209, \ldots & not in OEIS\\
000, 0100, 0101, 0111 & 10, 13, 18, 24, 30, 38, 49, 61, 75, 94, 117, 143,  \ldots & not in OEIS\\
0001, 0100, 0111, 100 & 10, 13, 18, 24, 32, 43, 57, 76, 101, 134, 178, 236, \ldots & not in OEIS\\
000, 0100, 0101, 1101 & 10, 13, 18, 25, 35, 48, 66, 91, 126, 174, 240, 331,  \ldots & not in OEIS\\
0011, 0101, 100 & 10, 13, 18, 25, 35, 50, 72, 104, 151, 220, 321, 469, \ldots & not in OEIS\\
0101, 1001, 110, 1101 & 10, 13, 18, 26, 37, 51, 70, 97, 135, 187, 258, 356,  \ldots & not in OEIS\\
000, 0101, 0110, 1101 & 10, 13, 19, 27, 37, 53, 74, 103, 146, 204, 286, 403, \ldots & not in OEIS\\
0001, 0101, 100 & 10, 13, 19, 28, 40, 58, 85, 124, 181, 265, 388, 568, \ldots & not in OEIS\\
000, 0100, 0101, 0110 & 10, 13, 19, 30, 47, 70, 102, 151, 228, 345, 517, 770, \ldots & not in OEIS\\
\hline
\end{tabular}
\caption{Forbidden words in the De Bruijn digraphs $B(2,4)$ with 10 vertices and 13 edges, and the sequence obtained with the numbers of vertices on $L^0(G)=G, L^1(G), L^2(G),\ldots$}
\label{table:forbid-words1}
\end{center}
\end{table}

\begin{table}[!ht]
\scriptsize
\begin{center}
\begin{tabular}{|c|l|l|}
\hline
Forbidden subwords & Sequence & OEIS  \\
\hline
0010, 011 & 11, 16, 22, 30, 41, 55, 74, 99, 132, 176, \dots & not in OEIS\\
0111, 100 & 11, 16, 23, 32, 44, 60, 81, 109, 146, 195,  \dots & not in OEIS\\
001, 0101, 100 & 11, 16, 23, 33, 47, 66, 93, 131, 183, 256, \dots & not in OEIS\\
000, 0100, 0111 & 11, 16, 24, 35, 49, 70, 100, 139, 195, 276, \dots & not in OEIS\\
000, 0110, 1001 & 11, 16, 24, 37, 56, 85, 128, 194, 293, 444, \dots & not in OEIS\\
000, 0101, 0110 & 11, 16, 25, 40, 63, 99, 155, 243, 382, 600,\dots & not in OEIS\\
000, 0011, 1011 & 11, 17, 24, 34, 47, 64, 87, 117, 157, 210,\dots & not in OEIS\\
0001, 011 & 11, 17, 25, 36, 51, 71, 98, 134, 182, 246,  \dots & not in OEIS\\
000, 0100, 1011 & 11, 17, 25, 36, 51, 72, 101, 141, 196, 272, \dots & not in OEIS\\
000, 1011, 1100 & 11, 17, 25, 37, 55, 82, 123, 185, 278, 418,\dots & not in OEIS\\
00, 0111, 1011 & 11, 17, 25, 38, 56, 83, 122, 180, 264, 388, \dots & not in OEIS\\
000, 1011, 1101 & 11, 17, 25, 39, 60, 92, 141, 216, 332, 509, \dots & not in OEIS\\
000, 0010, 1011 & 11, 17, 26, 39, 59, 89, 135, 204, 309, 467, \dots & not in OEIS\\
001, 0110 & 11, 17, 26, 40, 61, 93, 141, 214, 324, 491, \dots & not in OEIS\\
000, 0111, 1100 & 11, 17, 26, 40, 61, 93, 142, 216, 329, 501, \dots & not in OEIS\\
000, 0101, 1011 & 11, 17, 26, 40, 61, 94, 145, 223, 343, 528, \dots & not in OEIS\\
0000, 0101, 111 & 11, 17, 26, 41, 63, 97, 151, 234, 361, 559, \dots & not in OEIS\\
000, 0011, 1010 & 11, 17, 26, 41, 64, 99, 155, 242, 376, 587, \dots & not in OEIS\\
000, 0110, 1011 & 11, 17, 27, 42, 65, 101, 156, 242, 375, 581, \dots & not in OEIS\\
000, 0100, 0101 & 11, 17, 28, 46, 74, 119, 193, 313, 506, 818, \dots & not in OEIS\\
000, 0011, 0111 & 11, 18, 28, 43, 67, 102, 156, 239, 363, 554, \dots & not in OEIS\\
000, 0011, 1001 & 11, 18, 28, 46, 74, 120, 194, 314, 508, 822, \dots & not in OEIS\\
000, 0010, 1001 & 11, 18, 30, 48, 78, 126, 204, 330, 534, 864, \dots & not in OEIS\\
0001, 010 & 11, 18, 30, 49, 79, 128, 208, 337, 545, 882, \dots & not in OEIS\\
000, 0100, 100, 1010 & 11, 18, 30, 50, 83, 138, 229, 380, 631, \dots & not in OEIS\\
000, 0011, 1100 & 11, 18, 30, 50, 83, 138, 230, 383, 638, 1063 \dots & not in OEIS\\
000, 0011, 0100 & 11, 18, 30, 50, 84, 141, 236, 395, 661, \dots & not in OEIS\\
010, 1001 & 11, 19, 32, 53, 89, 149, 249, 417, 698, 1168 \dots & not in OEIS\\
\hline
\end{tabular}
\caption{Forbidden words in the De Bruijn digraphs $B(2,4)$ with 11 vertices and the sequence obtained with the numbers of vertices on $L^0(G)=G, L^1(G), L^2(G),\ldots$}
\label{table:forbid-words2}
\end{center}
\end{table}

\begin{table}[!ht]
\scriptsize
\begin{center}
\begin{tabular}{|c|l|l|}
\hline
Forbidden subwords & Sequence & OEIS  \\
\hline
0120, 0121, 0210 & 15, 19, 21, 25, 31, 38, 45, 55, \ldots & not in OEIS\\
0102, 0120, 0210 & 15, 19, 22, 27, 35, 43, 52, 65, \ldots & not in OEIS\\
0102, 0120, 0121 & 15, 20, 24, 29, 37, 44, 53, 65, \ldots & not in OEIS\\
0120, 0121, 0212 & 15, 20, 24, 29, 37, 45, 54, 66, \ldots & not in OEIS\\
0120, 0201, 0212 & 15, 20, 25, 31, 39, 46, 56, 69, \ldots & not in OEIS\\
0121, 0201, 0210 & 15, 20, 25, 31, 42, 56, 72, 94, \ldots & not in OEIS\\
0102, 0210, 0212 & 15, 20, 26, 33, 44, 57, 73, 96, \ldots & not in OEIS\\
0102, 0120, 0201 & 15, 20, 27, 37, 50, 67, 91, 124, \ldots & not in OEIS\\
0102, 0121, 0212 & 15, 21, 28, 35, 48, 63, 79, 108, \ldots & not in OEIS\\
020, 1021 & 15, 21, 28, 40, 55, 76, 104, 144, \ldots & not in OEIS\\
010, 0202, 0210 & 15, 21, 29, 41, 57, 80, 111, 155, \ldots & not in OEIS\\
0102, 0121, 0201 & 15, 21, 30, 41, 57, 81, 112, 155, \ldots & not in OEIS\\
020, 2101 & 15, 22, 31, 44, 62, 88, 125, 178,  \ldots & not in OEIS\\
0121, 212 & 15, 22, 31, 45, 64, 92, 132, 189,  \ldots & not in OEIS\\
1021, 212 & 15, 22, 31, 45, 65, 94, 135, 194,  \ldots & not in OEIS\\
1202, 212 & 15, 23, 34, 51, 76, 114, 170, 254, \ldots & not in OEIS\\
1201, 2102 & 16, 22, 28, 36, 46, 58, 72, 90, \ldots & not in OEIS\\
2012, 2102 & 16, 22, 28, 38, 52, 70, 92, 124, \ldots & not in OEIS\\
0120, 2120 & 16, 23, 31, 43, 60, 82, 112, 155, \ldots & not in OEIS\\
0120, 0212 & 16, 23, 31, 43, 60, 83, 114, 157, \ldots & not in OEIS\\
2101, 2120 & 16, 23, 31, 43, 61, 85, 118, 165, \ldots & not in OEIS\\
1021, 1210 & 16, 23, 32, 45, 63, 87, 121, 170, \ldots & not in OEIS\\
0102, 1201 & 16, 23, 32, 46, 67, 97, 139, 200, \ldots & not in OEIS\\
0210, 1021 & 16, 23, 33, 48, 68, 96, 137, 196, \ldots & not in OEIS\\
1202, 2010 & 16, 24, 34, 48, 68, 96, 136, 194, \ldots & not in OEIS\\
0102, 0121 & 16, 24, 34, 48, 69, 97, 137, 196, \ldots & not in OEIS\\
1020, 1202 & 16, 24, 34, 48, 70, 100, 142, 206, \ldots & not in OEIS\\
0201, 1202 & 16, 24, 34, 49, 70, 100, 144, 207, \ldots & not in OEIS\\
1201, 2010 & 16, 24, 34, 49, 71, 102, 146, 211, \ldots & not in OEIS\\
0121, 1020 & 16, 24, 34, 50, 74, 108, 158, 232, \ldots & not in OEIS\\
0102, 0212 & 16, 24, 35, 50, 74, 109, 158, 233, \ldots & not in OEIS\\
1012, 1210 & 16, 24, 36, 54, 80, 120, 180, 268, \ldots & not in OEIS\\
0212, 2021 & 16, 25, 36, 54, 81, 120, 180, 269, \ldots & not in OEIS\\
0201, 1020 & 16, 25, 38, 59, 90, 139, 214, 329, \ldots & not in OEIS\\
\hline
\end{tabular}
\caption{Forbidden words in the Square-free digraphs $SF(2,4)$ with 15 and 16 vertices and the sequence obtained with the numbers of vertices on $L^0(G)=G, L^1(G), L^2(G),\ldots$}
\label{table:forbid-words3}
\end{center}
\end{table}

\end{document}